\documentclass[reqno,a4wide,english]{amsart}

\usepackage{amsmath}

\usepackage[all]{xy}
\CompileMatrices
\newdir{ >}{{}*!/-10pt/\dir{>}}

\theoremstyle{plain}
\newtheorem{theorem}{Theorem}[section]

\theoremstyle{definition}

\newtheorem{emptythm}[theorem]{}

\newcommand{\OXx}{{\mathcal O}_{X,x}}

\newcommand{\OYy}{{\mathcal O}_{Y,y}}

\newcommand{\too}{\longrightarrow}

\newcommand{\ol}[1]{\overline{#1}}
\newcommand{\ul}[1]{\underline{#1}}

\DeclareMathOperator{\hgt}{ht}

\DeclareMathOperator{\Ker}{Ker}

\DeclareMathOperator{\id}{id}

\DeclareMathOperator{\Spec}{Spec}

\DeclareMathOperator{\rk}{\mathrm{rk}}

\DeclareMathOperator{\tensor}{\mathrm{T}}

\DeclareMathOperator{\PfM}{\mathrm{PfM}}

\DeclareMathOperator{\Kt}{\mathrm{K}}  

\newcommand{\MK}{\Kt^{M}}              

\newcommand{\MWK}{\Kt^{MW}}            
\newcommand{\MW}{\Kt^{W}}              
\newcommand{\WI}{\ul{W}}

\newcommand{\grFdI}{\ul{\iota}}

\DeclareMathOperator{\W}{\mathrm{W}}
\DeclareMathOperator{\GW}{\mathrm{GW}}
\DeclareMathOperator{\FdI}{\mathrm{I}}

\newcommand{\FdIunr}{\FdI_{unr}}

\newcommand{\Wunr}{\W_{unr}}

\DeclareMathOperator{\Orth}{\mathrm{O}}
\DeclareMathOperator{\ValGr}{\mathrm{D}}

\newcommand{\Z}{\mathbb{Z}}
\newcommand{\N}{\mathbb{N}}
\renewcommand{\P}{\mathbb{P}}
\newcommand{\A}{\mathbb{A}}

\newcommand{\Hp}{\mathbb{H}}

\newcommand{\maxid}{\mathfrak{m}}

\newcommand{\eg}{\textsl{e.g.}\ }
\newcommand{\ie}{\textsl{i.e.}\ }

\begin{document}

\title[Milnor-Witt $K$-groups of local rings]
{Milnor-Witt $K$-groups of local rings}

\author{Stefan Gille}
\email{gille@ualberta.ca}
\author{Stephen Scully}
\email{sscully@ualberta.ca}
\author{Changlong Zhong}
\email{zhongusc@gmail.com}
\address{Department of Mathematical and Statistical Sciences,
University of Alberta, Edmonton T6G 2G1, Canada}

\thanks{This work has been supported by an NSERC grant of the first author. The second
and third author are supported by PIMS postdoctoral fellowships.}


\date{February 26, 2015}

\begin{abstract}
We introduce Milnor-Witt $K$-groups of local rings and show that the $n$th
Milnor-Witt $K$-group of a local ring~$R$ which contains an infinite field of
characteristic not~$2$ is the pull-back of the $n$th power of the fundamental
ideal in the Witt ring of~$R$ and the $n$th Milnor $K$-group of~$R$ over the
$n$th Milnor $K$-group of~$R$ modulo~$2$. This generalizes the work of Morel-Hopkins
on Milnor-Witt $K$-groups of a field.
\end{abstract}

\maketitle

\section{Introduction}
\label{IntroSect}\bigbreak

\noindent
The Milnor-Witt $K$-theory of a field~$R$, denoted~$\MWK_{\ast}(R)=\bigoplus\limits_{n\in\Z}\MWK_{n}(R)$,
arises as an object of fundamental interest in motivic homotopy theory; namely, as the
``0-line'' part of the stable homotopy ring of the motivic sphere spectrum over~$R$, see
Morel~\cite{Mo04b}. Beginning from an initial presentation discovered in collaboration
with Hopkins, Morel~\cite{Mo04a} showed that, for a field~$R$, the group~$\MWK_{n}(R)$ 
is the pull-back of the diagram
\begin{equation}
\label{MWKpushoutEq}
\xymatrix{
 & \MK_{n}(R) \ar[d]^{e_{n}}
\\
\FdI^{n}(R) \ar[r] & \FdI^{n}(R)/\FdI^{n+1}(R)\, ,
}
\end{equation}
where $\MK_{n}(R)$ denotes the $n$th Milnor $K$-group of~$R$,
the symbol $\FdI^{n}(R)$ the $n$th power of fundamental ideal~$\FdI (R)$
in the Witt ring~$\W (R)$, and $e_{n}$ maps the
symbol~$\ell (a_{1})\cdot\ldots\cdot\ell (a_{n})$ to the class of the Pfister
form~$\ll a_{1},\ldots ,a_{n}\gg$. (Here~$\FdI^{n}(R)=\W (R)$ is understood
for~$n<0$.) In this form, Milnor-Witt $K$-groups of fields already appeared
implicitly in earlier work of Barge and Morel~\cite{BaMo99}, who used them
to introduce oriented Chow groups of an algebraic variety, a theory later
elaborated in the work of Fasel~\cite{Fa07,Fa08}. This so-called Chow-Witt theory was used
by Barge and Morel~\cite{BaMo99,BaMo00} to construct an Euler class invariant
for algebraic oriented vector bundles which, in analogy with its topological
counterpart, is the primary obstruction to the existence of a nowhere-vanishing
section, see Morel~\cite[Chap.\ 8]{A1AlgTop}. This circle of ideas has since been greatly
elaborated upon, see for instance the article~\cite{AsFa14} by Asok and Fasel for some
very recent developments.

\smallbreak

A further area where Milnor-Witt $K$-groups of a field show up is the theory
of framed motives, see for instance the recent work of Neshitov~\cite{Ne14}.

\medbreak

We introduce in this work Milnor-Witt $K$-groups~$\MWK_{\ast}(R)$
of a local ring~$R$. Our definition is the naive generalization of the Morel-Hopkins presentation given
in Morel's book~\cite[Def.\ 3.1]{A1AlgTop}, or Morel~\cite[Def.\ 5.1]{Mo04a}, \ie~$\MWK_{\ast}(R)$ is a
$\Z$-graded $\Z$-algebra generated by an element~$\hat{\eta}$ in degree~$-1$ and elements~$\{ a\}$ ($a$~a
unit in~$R$) in degree~$1$ modulo four relations, see Definition~\ref{MWKDef}. This definition has also
been considered in the recent work of Schlichting~\cite{Schl15}, where (extending earlier work of
Barge and Morel~\cite{BaMo00} and Hutchinson and Tao~\cite{HuTa10} over fields), it is shown that
Milnor-Witt $K$-groups arise as obstructions to integral homology stability for special linear groups
over local rings with infinite residue fields. This in turn is used then by Schlichting~\cite[\S 6]{Schl15}
to extend results of Morel on the vanishing of Euler class for oriented vector bundles over affine
schemes.

\medbreak

Our main result about these groups is the following, see~ Theorem~\ref{MWK-mainThm}.

\medbreak

\noindent
{\bfseries Theorem.}
{\it
Let~$R$ be a local ring which contains an infinite field of characteristic~$\not= 2$.
Then the $n$th Milnor-Witt $K$-group~$\MWK_{n}(R)$ is the pull-back of the
diagram~(\ref{MWKpushoutEq}) for all~$n\in\Z$.
}

\medbreak

\noindent
The proof uses a presentation of~$\FdI^{n}(R)$ given in~\cite{Gi15}, whose proof is based on
a recent deep theorem of Panin and Pimenov~\cite{PaPi10} on the existence of unimodular
isotropic vectors for quadratic forms over regular local rings.

\smallbreak

Another ingredient of our proof which is of some interest on its own is the following.
In Section~\ref{Chain-pSect} we show that over an arbitrary local ring, which contains~$\frac{1}{2}$,
two Pfister forms are isometric if and only if they are chain $p$-equivalent, see~\ref{Chain-pDef}
for the definition. This has been shown for fields by Elman and Lam~\cite{ElLa72} more than 40
years ago and seems to be new for local rings (although most likely known to experts). 

\smallbreak

In the last section we show for a regular local ring~$R$ containing an infinite field of
characteristic not~$2$ with quotient field~$K$ that the kernel of
$$
\MWK_{n}(K)\,\xrightarrow{\;(\partial_{P})_{\hgt P=1}\;}\,
\bigoplus\limits_{\hgt P=1}\MWK_{n-1}(R_{P}/PR_{P})\, ,
$$
where $\partial_{P}$ are the residue maps introduced by Morel~\cite[Sect.\ 3.2]{A1AlgTop},
is naturally isomorphic to~$\MWK_{n}(R)$ for all~$n\in\Z$. This implies via an argument
of Colliot-Th\'el\`ene~\cite[Sect.\ 2]{CT95} that the unramified Milnor-Witt
$K$-groups are birational invariants of smooth and proper schemes over an infinite field
of characteristic not~$2$.

\bigbreak

We want to point out that our proof of the main theorem does
not use results from Morel's work~\cite{Mo04a}, although we have
certainly borrowed ideas from there. Actually we prove the theorem
simultaneously also for fields of characteristic not~$2$.

\bigbreak

\noindent
{\bfseries Acknowledgement:}
The first author would like to thank Jean-Louis Colliot-Th\'el\`ene for showing him during
a conference in Mainz in September 2014 his proof that unramified Witt groups
are a birational invariant of smooth and proper schemes over a field. This was the starting
point for this work. We would further like to thank Jean Fasel, Detlev Hoffman, and Marco
Schlichting for very useful remarks and comments.

\bigbreak\bigbreak

\goodbreak
\section{Quadratic forms over local rings}
\label{QuadFormSect}\bigbreak

\begin{emptythm}
\label{ConventionSubSect}
We recall in this section some definitions and results of the algebraic theory of
quadratic forms over local rings. We refer for proofs and more information to
Scharlau~\cite[Chap.\ I, \S 6]{QHF}.

\smallbreak

We start with the following important

\smallbreak

\noindent
{\bfseries Convention.}
{\it Throughout this work we assume that all rings are commutative with~$1$ and
contain~$\frac{1}{2}$. In particular, we assume that fields are of characteristic
not~$2$.}
\end{emptythm}

\begin{emptythm}
\label{DefSubSect}
{\it Definitions.}
Let~$R$ be a local ring. A {\it quadratic form} over~$R$ is a map $q:V\too R$, where~$V$ is
a free $R$-module of finite rank, such that $q(\lambda v)=\lambda^{2}q(v)$ for all $\lambda\in R$
and~$v\in V$, and $b_{q}(v,w):=q(v+w)-q(v)-q(w)$ is a symmetric bilinear form on~$V$.
Throughout this work we assume as part of the definition that quadratic forms are non singular
(also called regular), \ie the associated bilinear form~$b_{q}$ is not degenerate. The
pair~$(V,q)$ is called a quadratic space.

\smallbreak

We denote for a quadratic space~$(V,q)$ the set of unit values of~$q$ by~$\ValGr (q)^{\times}$, \ie
$\ValGr (q)^{\times}:=\big\{ q(v)| v\in V\big\}\cap R^{\times}$, where~$R^{\times}$ is the multiplicative
group of units of~$R$.

\smallbreak

Two quadratic spaces $(V_{1},q_{1})$ and~$(V_{2},q_{2})$ are called {\it isomorphic} or {\it isometric}
if there exists an $R$-linear isomorphism $f:V_{1}\too V_{2}$, such that
$q_{2}(f(v))=q_{1}(v)$ for all~$v\in V_{1}$. Such a map~$f$ is called an isometry, and
we use the notation $(V_{1},q_{1})\simeq (V_{2},q_{2})$, or more briefly~$q_{1}\simeq q_{2}$,
to indicate that~$(V_{1},q_{1})$ and~$(V_{2},q_{2})$ are isometric. The group of all
automorphisms of the quadratic space~$(V,q)$ is denoted by $\Orth (V,q)$, or~$\Orth (q)$
only, and called the {\it orthogonal group} of~$(V,q)$ respectively of~$q$.

\smallbreak

As usual we denote the {\it orthogonal sum} and the {\it tensor product}
of two quadratic spaces~$(V_{1},q_{1})$ and~$(V_{2},q_{2})$ by
$(V_{1},q_{1})\perp (V_{2},q_{2})$ and $(V_{1},q_{1})\otimes (V_{2},q_{2})$,
respectively, or more briefly by~$q_{1}\perp q_{2}$ and~$q_{1}\otimes q_{2}$ only.
Note that the tensor product is defined since we assume that $\frac{1}{2}\in R$.

\smallbreak

Since~$\frac{1}{2}\in R$ every quadratic space has an orthogonal basis and
is therefore isomorphic to a diagonal form $<a_{1},\ldots ,a_{n}>$ for
appropriate~$a_{i}\in R^{\times}$. Another important consequence of our
assumption that~$2$ is invertible in~$R$ is {\it Witt cancellation}, \ie if
$q\perp q_{1}\simeq q\perp q_{2}$ for quadratic forms $q,q_{1},q_{2}$ over~$R$
then $q_{1}\simeq q_{2}$.
\end{emptythm}

\begin{emptythm}
\label{isotropicSubSect}
{\it Isotropic vectors.}
Let~$(V,q)$ be a quadratic space over the local ring~$R$. A vector~$v\in V\setminus\{ 0\}$,
such that $q(v)=0$, is called an {\it isotropic-} respectively if~$v$ is moreover unimodular
a {\it strictly isotropic vector}. If~$v$ is strictly
isotropic then there exists another isotropic vector~$w\in V$, such that~$b_{q}(v,w)=1$,
\ie~$v,w$ are a {\it hyperbolic pair} and $q|_{Rv\oplus Rw}\simeq <1,-1>$ is isometric to the
{\it hyperbolic plane}~$\Hp$. Then $(V,q)\simeq (V_{1},q|_{V_{1}})\perp\Hp$ for some
subspace~$V_{1}\subseteq V$.

\smallbreak

The quadratic space~$(V,q)$ respectively the quadratic form~$q$ is
called {\it isotropic} if there exists a strictly isotropic vector
(for~$q$) in~$V$. Otherwise~$(V,q)$ respectively~$q$ is called
{\it anisotropic}.
\end{emptythm}

\begin{emptythm}
\label{GrWittSubSect}
{\it Grothendieck-Witt groups.}
Let~$R$ be a local ring.
The Grothendieck group of isomorphism classes of quadratic spaces
over~$R$ with the orthogonal sum as addition is called
the {\it Grothendieck-Witt group} or {\it ring} of~$R$. It is in fact a commutative ring, where
the multiplication is induced by the tensor product. We denote this ring by~$\GW (R)$.
The quotient of~$\GW (R)$ by the ideal generated by the hyperbolic plane~$<1,-1>$ is
the {\it Witt group} or {\it Witt ring}~$\W (R)$ of~$R$. (Note that
since $\frac{1}{2}\in R$ every hyperbolic space over~$R$ is an orthogonal sum
of hyperbolic planes.)
\end{emptythm}

\begin{emptythm}
\label{FundamentalIdSubSect}
{\it The fundamental ideal of the Witt ring.}
Let~$R$ be a local ring. The {\it fundamental ideal}~$\FdI (R)$ of the Witt ring of~$R$
is the ideal consisting of the classes of even dimensional forms. It is additively
generated by the classes of $1$-Pfister forms~$\ll a\gg:=<1,-a>$. We denote its powers
by~$\FdI^{n}(R)$, $n\in\Z$, where $\FdI^{n}(R)=\W(R)$ for~$n\leq 0$ is understood.

\smallbreak

Obviously~$\FdI^{n}(R)$ is additively generated by $n$-Pfister forms
$$
\ll a_{1},\ldots ,a_{n}\gg\, :=\;\bigotimes\limits_{i=1}^{n}\ll a_{i}\gg\, ,
$$
$a_{1},\ldots ,a_{n}\in R^{\times}$, for all~$n\geq 1$. Note that a Pfister form~$q$
has an orthogonal decomposition $q=<1>\perp q'$. By Witt cancellation
the form~$q'$ is unique up to isometry and called the {\it pure subform} of the Pfister
form~$q$. We will use throughout the notation~$q'$ for the pure subform of a Pfister form~$q$.

\smallbreak

For later use we recall the following identities of Pfister forms.
\end{emptythm}

\begin{emptythm}
\label{isometryExpl}
{\bfseries Examples.}
Let~$R$ be a local ring.
\begin{itemize}
\item[(i)]
Let~$q=\ll a,b\gg$ be a $2$-Pfister form over~$R$. If $-c\in D(q')^{\times}$ then
we have $q\simeq\ll c,d\gg$ for some~$d\in R^{\times}$. This follows by comparing
determinants.

\smallbreak

\item[(ii)]
Assume that the two dimensional quadratic space $\ll a\gg =<1,-a>$ represents~$c\in R^{\times}$.
Then $<-b,ab>$ represents~$-bc$ for all $b\in R^{\times}$. Comparison of determinants then
implies $<-b,ab>\simeq <-bc,abc>$, and so we have then $\ll a,b\gg\simeq\ll a,bc\gg$.
\end{itemize}
\end{emptythm}

\begin{emptythm}
\label{GW-Z/2LemSubSect}
{\it A pull-back diagram.}
Another way to define the fundamental ideal is as follows. The {\it rank}
of a quadratic space~$(V,q)$ is by definition the rank of the underlying free
$R$-module~$V$. As the rank is additive on orthogonal
sums it induces a homomorphism $\rk:\GW (R)\too\Z$. Since hyperbolic spaces
have even rank this function induces in turn a homomorphism $\bar{\rk}:\W (R)\too\Z/2\Z$.
The fundamental ideal~$\FdI (R)$ is then the kernel of~$\bar{\rk}$.

\medbreak

As Witt cancellation holds for quadratic forms over~$R$ every quadratic space~$(V,q)$
decomposes up to isometry uniquely $(V,q)\simeq (W,\varphi)\perp (U,\phi)$ with~$\varphi$
anisotropic and~$(U,\phi)$ hyperbolic. This implies the following well known fact.

\medbreak

\noindent
{\bfseries Lemma.}
{\it
Let~$R$ be a local ring. Then the diagram
$$
\xymatrix{
\GW (R) \ar[r]^-{\rk} \ar[d] & {\Z} \ar[d]
\\
\W (R) \ar[r]_-{\bar{\rk}} & {\Z/2\Z}
}
$$
is a pull-back diagram.
}
\end{emptythm}

\begin{emptythm}
\label{ReflcationSubSect}
{\it Reflections.}
A vector~$0\not= v\in V$ which is not isotropic is called {\it anisotropic},
respectively {\it strictly anisotropic} if~$q(v)$ is a unit in~$R$. A strictly anisotropic
vector has to be unimodular. If~$v\in V$ is a strictly anisotropic vector for the
quadratic form~$q$ then
$$
\tau_{v}\, :\; V\,\too\, V\, ,\; x\,\longmapsto\, x-\frac{b_{q}(v,x)}{q(v)}\cdot v
$$
is a well defined $R$-linear map, called the {\it reflection} associated with
the vector~$v$. Recall that we have $q(\tau_{v}(x))=q(x)$ for all~$x\in V$,
\ie~$\tau_{v}$ is an isometry and so defines an element of~$\Orth (q)$.
\end{emptythm}

\begin{emptythm}
\label{ReductionSubSect}
{\it Reduction.}
We continue with the notation of the last section.
Given a quadratic space~$(V,q)$ over~$R$ then~$(\bar{V},\bar{q}):=k\otimes_{R}(V,q)$
is a quadratic space over the residue field~$k$ of~$R$, the {\it canonical reduction}
of~$(V,q)$. Note that if~$(V,q)$ is isotropic then also the reduction~$(\bar{V},\bar{q})$
is isotropic.

\smallbreak

If~$f:\phi\simeq\psi$ is an isometry then $\bar{f}:=\id_{k}\otimes f$
is an isometry between~$\bar{\phi}$ and~$\bar{\psi}$. In particular, we have then
a homomorphism of orthogonal groups
$$
\rho_{\phi}\, :\;\Orth (\phi)\,\too\,\Orth (\bar{\phi})\, ,\; \alpha\,\longmapsto\,\bar{\alpha}\, .
$$
This map is surjective. In fact, by the Cartan-Dieudonn\'e Theorem, see \eg~\cite[Chap.\ 1, Thm.\ 5.4]{QHF},
the group~$\Orth (\bar{q})$ is generated by reflections~$\tau_{\bar{v}}$ with~$\bar{v}\in\bar{V}$
an anisotropic vector. If~$v\in V$ is a vector which maps to~$\bar{v}$ under the quotient
map~$V\too\bar{V}$ then~$v$ is a strictly anisotropic vector in~$V$, and therefore
the reflection~$\tau_{v}\in\Orth (q)$ exists. We have~$\bar{\tau}_{v}=\rho_{\phi}(\tau_{v})=\tau_{\bar{v}}$
which proves the claim.
\end{emptythm}

\goodbreak
\section{Chain $p$-equivalent and isometric quadratic forms over local rings}
\label{Chain-pSect}\bigbreak

\begin{emptythm}
\label{Chain-pDef}
{\it Definition of chain $p$-equivalence.}
Chain $p$-equivalence of quadratic forms over fields has been introduced by
Elman and Lam~\cite{ElLa72}. The definition carries over to rings word by word.

\smallbreak

\noindent
{\bfseries Definition.}
Let~$R$ be a local ring with maximal ideal~$\maxid$ and residue field
$k=R/\maxid$. Let~$n\geq 2$ be an integer, and $\phi=\ll a_{1},\ldots ,a_{n}\gg$
and $\varphi=\ll b_{1},\ldots ,b_{n}\gg$ two $n$-Pfister forms.
\begin{itemize}
\item[(i)]
The $n$-Pfister forms~$\phi$ and~$\varphi$ are called {\it simply $p$-equivalent}
if there exist indices $1\leq i< j\leq n$, such that $\ll a_{i},a_{j}\gg\simeq\ll b_{i},b_{j}\gg$
and $a_{l}=b_{l}$ for all $l\not= i,j$.

\smallbreak

\item[(ii)]
The $n$-Pfister forms~$\phi$ and~$\varphi$ are called {\it chain $p$-equivalent} if
there exists a chain $\phi=\mu_{0},\mu_{1},\ldots ,\mu_{r}=\varphi$ of $n$-Pfister forms
over~$R$, such that $\mu_{i}$ is simply $p$-equivalent to~$\mu_{i+1}$ for all
$0\leq i\leq r-1$.
\end{itemize}

\noindent
Following Elman and Lam~\cite{ElLa72} we use the notation $\phi\approx\varphi$
to indicate that~$\phi$ and~$\varphi$ are chain $p$-equivalent.

\smallbreak

The aim of this section is to show that over a local ring~$R$, where~$2$ is a unit,
Pfister forms are isometric if and only if they are chain $p$-equivalent. In the case
where~$R$ is a field of characteristic not~$2$, this was first proved by Elman and
Lam~\cite{ElLa72}. To treat the more general case, we follow essentially their
arguments, but some modifications are necessary, mainly for the reason that
a non-zero element in a general local ring need not to be a unit. We overcome
these obstructions using Lemma~\ref{maintoolLem}
respectively Lemma~\ref{3-5Lem} if the residue field is small.
\end{emptythm}

\begin{emptythm}
\label{maintoolLem}
{\bfseries Lemma.}
{\it
Let~$R$ be a local ring whose residue field~$k$ is not the field with~$3$ or~$5$ elements.
Let~$(V,\phi)$ and~$(W,\varphi)$ be quadratic spaces over~$R$, and set $q:=\phi\perp\varphi$.
Then given~$a\in\ValGr (q)^{\times}$ there exists~$v\in V$ and~$w\in W$, such that
\begin{itemize}
\item[(a)]
$a\, =\, q(v,w)\, =\,\phi (v)+\varphi (w)$, and

\smallbreak

\item[(b)]
both~$\phi (v)$ and~$\varphi (w)$ are units in~$R$.
\end{itemize}
}

\begin{proof}
Let~$x=(v',w')$. We assume first that~$R=k$ is a field.

\smallbreak

If~$\phi (v')\not= 0\not=\varphi (w')$
there is nothing to prove, so assume one of these values is zero, say $\varphi (w')=0$.
Then $a=\phi (v')=q(v')$, and so we can assume that $x=v'$.

\smallbreak

The quadratic form~$\varphi$ is non singular and so there exists~$z\in W$,
such that $\varphi (z)\not= 0$. As~$R=k$ has at least~$7$ elements
there exists~$\lambda_{0}\in k^{\times}$,
such that
$$
\lambda^{2}_{0}\,\not=\,\pm\frac{\phi (v')}{\varphi (z)}\, .
$$
Consider the vector $u:=v'+\lambda_{0}\cdot z$. We have
$q(u)=\phi (v')+\lambda^{2}_{0}\cdot\varphi(z)$ which is not~$0$ by our
choice of~$\lambda_{0}$, \ie~$u$ is an anisotropic vector and so the
reflection~$\tau_{u}\in\Orth (q)$ is defined. Set now
$$
(v,w)\, :=\;\tau_{u}(v')\, .
$$
As~$\tau_{u}$ is in~$\Orth (q)$ we have
$a=q(v')=q(\tau_{u}(v'))=q(v,w)=\phi (v)+\varphi (w)$, and so~(a).
We are left to show that~$\phi(v)$ and~$\varphi (w)$ are both non zero. For this
we compute
$$
(v,w)\, =\,\tau_{u}(v')\, =\,
\left(1-\frac{2\phi (v')}{\phi (v')+\lambda_{0}^{2}\varphi (z)}\right)\cdot v'\, +\,
\frac{2\lambda_{0}\phi (v')}{\phi (v')+\lambda_{0}^{2}\varphi (z)}\cdot z\, ,
$$
and so $v=\left(1-\frac{2\phi (v')}{\phi (v')+\lambda_{0}^{2}\varphi (z)}\right)\cdot v'
=:c\cdot v'$ and $w=\frac{2\lambda_{0}\phi (v')}{\phi (v')+\lambda_{0}^{2}\varphi (z)}
\cdot z=:d\cdot z$. By our choice of~$\lambda_{0}$ both coefficients~$c,d$ are non zero and so
$$
\phi (v)=c^{2}\cdot\phi (v')\,\not=\, 0\,\not=\,
d^{2}\varphi (z)=\varphi (w)
$$
as desired.

\medbreak

We come now to the general case, \ie~$R$ is a local ring whose residue field~$k$
has at least $7$~elements. If~$U$ is an $R$-module we denote by~$\bar{u}$ the
image of~$u\in U$ in $\bar{U}=k\otimes_{R}U$.

\smallbreak

Since~$a\in R^{\times}$ its residue~$\bar{a}$ in~$k$ is non zero. By the
field case there are then~$\bar{v}\in\bar{V}$ and~$\bar{w}\in\bar{W}$, such that
$\bar{a}=\bar{q}(\bar{v},\bar{w})=\bar{\phi}(\bar{v})+\bar{\varphi}(\bar{w})$, and
$$
\bar{\phi}(\bar{v})\,\not= 0\,\not=\,\bar{\varphi}(\bar{w})\, .
$$
Since $0\not=\bar{a}=\bar{q}(\bar{v},\bar{w})=\bar{q}(\bar{v}',\bar{w}')$ there
exists~$\bar{\tau}\in\Orth (\bar{q})$, such that
$$
\bar{\tau}(\bar{v}',\bar{w}')\, =\, (\bar{v},\bar{w})\, .
$$
As seen in~\ref{ReductionSubSect} there exists $\tau\in\Orth (q)$ whose canonical
reduction is~$\bar{\tau}$. We set then $(v,w):=\tau (v',w')$. Clearly we have
then $a=q(v,w)$, and~$\phi (v)$ and~$\varphi (w)$ are both units in~$R$ since
their reductions $\ol{\phi (v)}=\bar{\phi}(\bar{v})$ and~$\ol{\varphi (w)}=
\bar{\varphi}(\bar{w})$ are both non zero in the residue field~$k$.
\end{proof}

\smallbreak

\noindent
{\bfseries Remark.}
The form $q=<1>\perp <1>$ and $a=1$ shows that the lemma does not hold
for~$R$ a field with~$3$ or~$5$ elements.

\smallbreak

To handle also the case of local rings whose residue fields have only~$3$ or~$5$
elements we prove a more specialized lemma.
\end{emptythm}

\begin{emptythm}
\label{3-5Lem}
{\bfseries Lemma.}
{\it
Let~$R$ be a local ring whose residue field~$k$ has~$3$ or~$5$ elements, and~$(W,\varphi)$
a quadratic space of rank~$\geq 3$ over~$R$.
\begin{itemize}
\item[(i)]
Let $q=bt^{2}+\varphi$ with~$b\in R^{\times}$, and~$a\in\ValGr (q)^{\times}$.
Then there is $s\in R$ and~$w\in W$, such that $\varphi (w)\in R^{\times}$
and $a=bs^{2}+\varphi (w)$ (note that we do not claim that~$s$ is a unit in~$R$).

\smallbreak

\item[(ii)]
Let~$(V,\phi)$ be another quadratic space of rank~$\geq 3$ over~$R$ and
$q=\phi\perp\varphi$. Let $a\in\ValGr (q)^{\times}$. Then there are vectors~$v\in V$
and~$w\in W$, such that both~$\phi (v)$ and~$\varphi (w)$ are units in~$R$ and
$a=\phi (v)+\varphi (w)$. 
\end{itemize}
}

\begin{proof}
The proof uses the fact that over a finite field of characteristic not~$2$ every quadratic
form of rank~$\geq 3$ is isotropic, see \eg~\cite[Chap.\ 2, Thm.\ 3.8]{QHF}, and
so is universal, \ie represents every element of the field.

\smallbreak

For~(i) let $t_{0}\in R$ and~$w_{0}\in V$, such that $a=bt_{0}^{2}+\varphi (w_{0})$.
Then $0\not=\bar{a}=\bar{b}\bar{t}_{0}+\bar{\varphi}(\bar{w}_{0})$. Since $\dim\bar{\varphi}\geq 3$
there exists by the remark above $\bar{w}\in\bar{W}$, such that $\bar{\varphi}(\bar{w})=\bar{a}$.
Then there is~$\bar{\tau}\in\Orth (\bar{q})$, such that $\tau (\bar{t_{0}},\bar{w}_{0})=(0,\bar{w})$.
Let~$\tau\in\Orth (q)$ be a preimage of~$\bar{\tau}$ under the reduction map~$\Orth (q)\too\Orth (\bar{q})$,
see~\ref{ReductionSubSect}. The vector $(s,w):=\tau (t_{0},w_{0})$ does the job.

\smallbreak

For~(ii) we use a similar reasoning. There are~$v_{0}\in V$ and~$w_{0}\in W$, such that
$a=\phi (v_{0})+\varphi (w_{0})$. Since~$k$ has at least three elements and~$\bar{\phi}$
is universal by dimension reasons, there exists~$\bar{v}\in\bar{V}$, such that
$\bar{\phi}(\bar{v})\not= 0\not= \bar{a}-\bar{\phi}(\bar{v})$. By the universality of~$\bar{\varphi}$
there exists then~$\bar{w}\in\bar{W}$, such that $\bar{\varphi}(\bar{w})=\bar{a}-\bar{\phi}(\bar{v})$.
Let then~$\bar{\tau}$ be an automorphism of~$q$, such that $\bar{\tau}(\bar{v}_{0},\bar{w}_{0})=
(\bar{v},\bar{w})$, and~$\tau\in\Orth (q)$ a preimage of~$\bar{\tau}$ under the reduction
map. The vector~$(v,w):=\tau (v_{0},w_{0})$ has the desired properties.
\end{proof}

We begin now the proof of the main result of this section. Except for the use of
Lemmas~\ref{maintoolLem} and~\ref{3-5Lem} above the arguments are almost word by
word the same as in the field case (in fact even a bit shorter since these lemmas
allow us to avoid some case by case considerations).
For the sake of completeness (and to convince the reader of the
correctness of the results) we give the details.

\smallbreak

We start with the following lemma which corresponds to Elman and
Lam~\cite[Prop.\ 2.2]{ElLa72}.
\end{emptythm}

\begin{emptythm}
\label{ValueLem}
{\bfseries Lemma.}
{\it
Let~$q=\ll a_{1},\ldots ,a_{n}\gg$ be a $n$-Pfister form over~$R$, $n\geq 2$.
If~$-b\in\ValGr (q')^{\times}$, where~$q'$ is the pure subform of~$q$,
then there exist $b_{2},\ldots ,b_{n}\in R^{\times}$, such that
$q\,\approx\,\ll b,b_{2},\ldots ,b_{n}\gg$.
}

\begin{proof}
We prove this by induction on~$n\geq 2$. The case~$n=2$ is Example~\ref{isometryExpl}~(i),
so let~$n\geq 3$. Set $\phi=\ll a_{1},\ldots ,a_{n-1}\gg$. Then by Witt cancellation
we have $q'\simeq\phi'\perp <-a_{n}>\otimes\,\phi$.
By Lemmas~\ref{maintoolLem} and~\ref{3-5Lem} we can assume that
$-b=-x-a_{n}y$ and~$y=t^{2}-z$ with $-x,-z\in\ValGr (\phi')^{\times}$ and~$y\in\ValGr (\phi)^{\times}$.
By induction we have then
$$
\ll x,b_{2},\ldots ,b_{n-1}\gg\,\approx\,\phi\,\approx\,
\ll z,c_{2},\ldots ,c_{n-1}\gg
$$
for some $b_{i},c_{i}\in R^{\times}$, $2\leq i\leq n-1$. As~$<1,-z>$ represents~$y$ we
conclude then by Example~\ref{isometryExpl}~(ii) that $\ll z,a_{n}\gg\simeq\ll z,a_{n}y\gg$
and therefore
$$
\begin{array}{r@{\,\approx\,}l}
q\, =\,\phi\,\otimes\ll a_{n}\gg  & \ll z,c_{2},\ldots ,c_{n-1},a_{n}\gg\,\approx\,
            \ll z,c_{2},\ldots ,c_{n-1},a_{n}y\gg \\[2mm]
 & \ll x,b_{2},\ldots ,b_{n-1},a_{n}y\gg\, .
\end{array}
$$
The pure subform of the $2$-Pfister form~$\ll x,a_{n}y\gg$ represents~$-b$
and so by the case~$n=2$ we have $\ll x,a_{n}y\gg\simeq\ll b,b_{n}\gg$
for some~$b_{n}\in R^{\times}$. We are done.
\end{proof}

The following two assertions correspond to~\cite[Cor.\ 2.5 and Thm.\ 2.6]{ElLa72}.
\end{emptythm}

\begin{emptythm}
\label{ValueCor-Lem}
{\bfseries Lemma.}
{\it
Let $p=\ll a_{1},\ldots ,a_{m}\gg$ and $q=\ll b_{1},\ldots ,b_{n}\gg$, $m,n\geq 1$ be
two Pfister forms over the local ring~$R$. Then:
\begin{itemize}
\item[(i)]
If $c\in\ValGr (p)^{\times}$ then for all $d\in R^{\times}$ we have
$p\,\otimes\ll d\gg\approx p\,\otimes\ll cd\gg$.

\smallbreak

\item[(ii)]
If~$-c\in\ValGr (p\,\otimes q')^{\times}$ then there are units
$c_{2},\ldots ,c_{n}\in R$, such that
$$
p\,\otimes q\,\approx\,p\,\otimes\ll c,c_{2},\ldots ,c_{n}\gg\, .
$$
\end{itemize}
}

\begin{proof}
We start with~(i). If~$m=1$ this follows form Example~\ref{isometryExpl}~(ii),
so let~$m\geq 2$. We can assume by Lemmas~\ref{maintoolLem} and~\ref{3-5Lem} that
$c=t^{2}-z$ with~$-z\in\ValGr (p')^{\times}$, and so by the lemma
above $p\,\otimes\ll d\gg$ is chain $p$-equivalent to $\ll z,c_{2},\ldots ,c_{m},d\gg$
for some~$c_{i}\in R^{\times}$, $2\leq i\leq m$. As $c\in\ValGr (<1,-z>)^{\times}$
we have by Examples~\ref{isometryExpl}~(ii) the isometry $\ll z,d\gg\simeq\ll z,cd\gg$
from which the claim follows.

\smallbreak

We prove now~(ii) by induction on $n\geq 1$. The case~$n=1$ is an immediate
consequence of~(i) as then $-c=-b_{1}\cdot x$ with~$x\in\ValGr (p)$. So
let $n\geq 2$. Write $q=\phi\,\otimes\ll b_{n}\gg$, and so
$p\,\otimes q'=p\,\otimes\phi'\perp <-b_{n}>\otimes\, p\,\otimes\,\phi$.
Again by Lemmas~\ref{maintoolLem} and~\ref{3-5Lem} we can then write
$-c=-x-b_{n}\cdot y$ with $-x\in\ValGr (p\,\otimes\phi')^{\times}$ and $y\in\ValGr (p\otimes\phi)^{\times}$.
By induction we have then $p\,\otimes\phi\approx p\,\otimes\ll x,c_{2},\ldots ,c_{n-1}\gg$
for some $c_{2},\ldots, c_{n-1}\in R^{\times}$. On the other hand since
$y\in\ValGr (p\,\otimes\phi)^{\times}$ we have by~(i) that
$p\,\otimes q=p\,\otimes\phi\,\otimes\ll b_{n}\gg\approx p\,\otimes\phi\,\otimes\ll b_{n}y\gg$, and so
$$
p\,\otimes q\approx p\,\otimes\ll x,c_{2},\ldots ,c_{n-1},b_{n}y\gg\, .
$$
From this the claim follows since the pure subform of the $2$-Pfister form
$\ll x,b_{n}y\gg$ represents the unit~$-c$, and so by Lemma~\ref{ValueLem}
above we have $\ll x,b_{n}y\gg\simeq\ll c,c_{n}\gg$ for some~$c_{n}\in R^{\times}$.
\end{proof}

We prove now the main result of this section.
\end{emptythm}

\begin{emptythm}
\label{p-chainThm}
{\bfseries Theorem.}
{\it
Let~$n\geq 2$, and~$p=\ll a_{1},\ldots ,a_{n}\gg$ and~$q=\ll b_{1},\ldots ,b_{n}\gg$
be two $n$-Pfister forms over the local ring~$R$. Then
$$
p\simeq q\quad\Longleftrightarrow\quad p\approx q\, .
$$
}

\begin{proof}
Obviously $p\approx q$ implies~$p\simeq q$. For the other direction
the case~$n=2$ is by definition, so let~$n\geq 3$. By Witt
cancellation the pure subforms of~$p$ and~$q$ are isomorphic,
and so $-b_{1}\in\ValGr (p')^{\times}$.

\smallbreak

Hence by Lemma~\ref{ValueLem} we have
$p\approx\ll b_{1},c_{2},\ldots ,c_{n}\gg$
for appropriate $c_{2},\ldots ,c_{n}\in R^{\times}$. Let $1\leq r\leq n$
be maximal, such that $p\approx\ll b_{1},\ldots ,b_{r},d_{r+1},\ldots ,d_{n}\gg$
for some units $d_{r+1},\ldots ,d_{n}$ in~$R$. We claim that $r=n$ which finishes the
proof. Assume the contrary, \ie~$r<n$. We set then
$$
\phi =\ll d_{r+1},\ldots ,d_{n}\gg\qquad\mbox{and}\qquad
\varphi =\ll b_{r+1},\ldots ,b_{n}\gg\, .
$$
Then we have
$$
\ll b_{1},\ldots ,b_{r}\gg\otimes\phi\,\simeq\,
\ll b_{1},\ldots ,b_{r}\gg\otimes\varphi\, ,
$$
and so by Witt cancellation we have $\ll b_{1},\ldots ,b_{r}\gg\otimes\,\phi'
\simeq\ll b_{1},\ldots ,b_{r}\gg\otimes\,\varphi'$. Therefore $-b_{r+1}\in\ R^{\times}$
is represented by $\ll b_{1},\ldots ,b_{r}\gg\otimes\,\phi'$ and so it follows
from Lemma~\ref{ValueCor-Lem}~(ii) that
$$
p\,\approx\,\ll b_{1},\ldots ,b_{r},d_{r+1},\ldots ,d_{n}\gg\,\approx\,
\ll b_{1},\ldots ,b_{r},b_{r+1},e_{r+2},\ldots ,e_{n}\gg
$$
for some $e_{r+2},\ldots ,e_{n}\in R^{\times}$, contradicting the maximality of~$r$.
We are done.
\end{proof}
\end{emptythm}

\goodbreak
\section{Witt $K$-theory of a local ring}
\label{WKSect}\bigbreak

\begin{emptythm}
\label{7-Convention}
{\bfseries Convention.}
{\it
In this section~$R$ denotes a local ring with residue field~$k$. If~$R$ is not
a field, \ie $R\not= k$, we assume that~$k$ has at least $5$~elements.
}
\end{emptythm}

\begin{emptythm}
\label{WKDefSubSect}
{\it Witt $K$-theory of a local ring.}
Witt $K$-theory of fields has been introduced by Morel~\cite{Mo04a}. Our definition for
a local ring is the obvious and straightforward generalization.

\smallbreak

\noindent
{\bfseries Definition.}
The {\it Witt $K$-ring of~$R$} is the quotient of the graded
and free $\Z$-algebra generated by elements~$[a]$ ($a\in R^{\times}$) in degree~$1$
and one element~$\eta$ in degree~$-1$ by the two sided ideal, which is generated by
the expressions
\begin{itemize}
\item[(WK1)]
$\eta\cdot [a]-[a]\cdot\eta$ with $a\in R^{\times}$;

\smallbreak

\item[(WK2)]
$[ab]-[a]-[b]+\eta\cdot [a]\cdot [b]$ with $a,b\in R^{\times}$;

\smallbreak

\item[(WK3)]
$[a]\cdot [1-a]$ with~$a\in R^{\times}$, such that~$1-a$ in~$R^{\times}$; and

\smallbreak

\item[(WK4)]
$2-\eta\cdot [-1]$.
\end{itemize}
We denote this graded $\Z$-algebra by $\MW_{\ast}(R)=\bigoplus\limits_{n\in\Z}\MW_{n}(R)$.
Note that by~(WK2) the $n$th graded piece~$\MW_{n}(R)$ is generated by all
products $[a_{1}]\cdot\ldots\cdot [a_{n}]$ for~$n\geq 1$.

\smallbreak

Following Morel~\cite{Mo04a} we set
$$
<a>_{\MW}\, :=\; 1-\eta\cdot [a]\in\MW_{0}(R)\, .
$$
A straightforward computation using~(WK1) and~(WK2) shows that
$$
<ab>_{\MW}\, =\, <a>_{\MW}\,\cdot\, <b>_{\MW}
$$
for all $a,b\in R^{\times}$. Note that using this notation~(WK2)
can be reformulated as
$$
[ab]\, =\, [a]+<a>_{\MW}\cdot [b]\, .
$$
\end{emptythm}

\begin{emptythm}
\label{elemIdSubSect}
{\it Some elementary identities in~$\MW_{\ast}(R)$.}
The following identities are proven for fields in Morel's article~\cite{Mo04a}.
For the convenience of our reader we recall the rather short arguments which
also work for local rings.
\begin{itemize}
\item[(1)]
By~(WK4) we have $<-1>_{\MW}=-1$, and so using
$$
<1>_{\MW}\, =\, <-1>_{\MW}\cdot <-1>_{\MW}
$$
we have~$<1>_{\MW}=1$. The later equation implies then~$[1]=0$
since by~(WK2) we have $[1]=[1\cdot 1]=[1]+<1>_{\MW}\cdot [1]$.

\smallbreak

It follows from this that~$<a>_{\MW}$ is invertible
with inverse
$$
<a>_{\MW}^{-1}\, =\, <a^{-1}>_{\MW}
$$
for all~$a\in R^{\times}$.

\smallbreak

\item[(2)]
By~(WK2) we have
$$
\eta\cdot [a]\cdot [b]\, =\,\eta\cdot [b]\cdot [a]
$$
for all~$a,b\in R^{\times}$. This implies in particular that
$$
<a>_{\MW}\cdot [b]\, =\, [b]\cdot <a>_{\MW}
$$
for all~$a,b\in R^{\times}$.

\smallbreak

\item[(3)]
We have
\begin{itemize}
\item[(a)]
$[b^{-1}]\, =\, -<b^{-1}>_{\MW}\cdot [b]$, and

\smallbreak

\item[(b)]
$[ab^{-1}]=[a]-<ab^{-1}>_{\MW}\cdot [b]$
\end{itemize}
for all $a,b\in R^{\times}$. In fact,
we have by~(1) above $[1]=0$, and so $0=[b^{-1}\cdot b]=[b^{-1}]+<b^{-1}>_{\MW}\cdot [b]$,
hence~(a).

\smallbreak

Since $<ab^{-1}>_{\MW}=<a>_{\MW}\cdot <b^{-1}>_{\MW}$ we get from this
$$
[a]-<ab^{-1}>_{\MW}\cdot [b]=[a]+<a>_{\MW}\cdot [b^{-1}]=[ab^{-1}]
$$
as claimed.
\end{itemize}
The proof of the following lemma uses our assumption that if~$R$ is not a
field then the residue field has at least $5$~elements.
\end{emptythm}

\begin{emptythm}
\label{a-aRelLem}
{\bfseries Lemma.}
{\it
We have $[-a]\cdot [a]=0$ for all $a\in R^{\times}$.
}

\begin{proof}
The argument below is an adaption of the one in Nesterenko and Suslin~\cite{NeSu90}
for Milnor $K$-theory of local rings.

\smallbreak

Assume first that~$a$, $1-a$ and~$1-a^{-1}$ are units in~$R$, \ie the
residue~$\bar{a}\in k$ is not~$1$. Then $-a=\frac{1-a}{1-a^{-1}}$, and so
we have
$$
[-a]=[1-a]-<-a>_{\MW}\cdot [1-a^{-1}]
$$
using~\ref{elemIdSubSect}~(3) above. Therefore since $[1-a][a]=0$ we have
$$
\begin{array}{r@{\, =\,}l}
   [-a]\cdot [a] & -<-a>_{\MW}\cdot [1-a^{-1}]\cdot [a] \\[2mm]
     & <-a>_{\MW}\cdot [1-a^{-1}]\cdot [a^{-1}]\cdot <a>_{\MW}\, =\, 0\, ,
\end{array}
$$
(using~\ref{elemIdSubSect}~(2) and~(3)), and analogous $[a]\cdot [-a]=0$.
This proves the result in particular if~$R$ is a field.

\smallbreak

Assume now that~$R$ is not a field and~$\bar{a}=1$ in the residue field~$k$.
By our assumption~$k$ has at least $5$~elements. Let~$c\in R^{\times}$,
such that $\bar{c}\not= 1$ in~$k$. Then also $\bar{a}\cdot\bar{c}\not= 1$ and so
as we have shown already $[-ac]\cdot [ac]=[-c]\cdot [c]=[c]\cdot [-c]=0$.
This implies using~(WK2) together with~\ref{elemIdSubSect}~(2)
\begin{equation}
\label{a-aEq}
[-a]\cdot [a]\, =\, -[c]\cdot [a]-[a]\cdot [c]+\eta\cdot
\big([-a]\cdot [a]+[a]\cdot [a]\big)\cdot [c]\, .
\end{equation}
Since $|k|\geq 5$ there exists~$e,f\in R^{\times}$, such that $\bar{e},\bar{f}\not= 1$
and also $\bar{e}\bar{f}\not= 1$. Applying Equation~(\ref{a-aEq}) to $c=ef$ and
using~\ref{elemIdSubSect}~(2) we get
$$
\begin{array}{r@{\,}c@{\,}l}
  [-a]\cdot [a] & = & -[ef]\cdot [a]-[a]\cdot [ef]+\eta\cdot
\big([-a]\cdot [a]+[a]\cdot [a]\big)\cdot [ef] \\[4mm]
 & = & -[e]\cdot [a]-[a]\cdot [e]+\eta
    \cdot\big([-a]\cdot [a]+[a]\cdot [a]\big)\cdot [e] \\[3mm]
 & & \;\;\;\; -[f]\cdot [a] -[a]\cdot [f]+\eta\cdot
      \big( [-a]\cdot [a]+[a]\cdot [a]\big)\cdot [f] \\[3mm]
 & & \;\;\;\;\;\;\;\;\;\; +2\eta\cdot [e]\cdot [f]\cdot [a]\, -\eta^{2}\cdot ([-a]\cdot [a]+[a]\cdot [a])\cdot [e]\cdot [f] \\[4mm]
 & = & [-a]\cdot [a]+[-a]\cdot [a]\, +
         \eta\cdot [e]\cdot [f]\cdot [a]\cdot\big(2-\eta\cdot [-a]+\eta\cdot [a]\big)\, ,
\end{array}
$$
where the last equation is by~(\ref{a-aEq}) for $c=e,f$ and~\ref{elemIdSubSect}~(2).
Now by~(WK2) and~\ref{elemIdSubSect}~(2) we have
$$
\begin{array}{r@{\; =\;}l}
   2-\eta\cdot [-a]-\eta\cdot [a] & 2-\eta\cdot [-1]-\eta\cdot [a]+\eta^{2}\cdot [-1]\cdot [a]-\eta\cdot [a] \\[4mm]
      & 2-\eta\cdot [-1]+\eta\cdot [a]\cdot (-2+\eta\cdot [-1])\, ,
\end{array}
$$
which is~$0$ by~(WK4). Hence $[-a]\cdot [a]=0$. We are done.
\end{proof}
\end{emptythm}

\begin{emptythm}
\label{a-1RelCor}
{\bfseries Corollary.}
{\it
We have
\begin{itemize}
\item[(i)]
$[a]\cdot [a]=[a]\cdot [-1]$ for all~$a\in R^{\times}$; and

\smallbreak

\item[(ii)]
$[ab^{2}]=[a]$, and so in particular $[b^{2}]=0$ for all $a,b\in R^{\times}$.
\end{itemize}
}

\begin{proof}
This can be proven as in the field case, see~\cite{Mo04a}. For
the sake of completeness we recall briefly the details.

\smallbreak

For~(i), we have $[-a]=[-1]+<-1>_{\MW}\cdot [a]=[-1]-[a]$
by~(WK2) and~\ref{elemIdSubSect}~(1), and so using the lemma above
$0=[a]\cdot [-a]=[a]\cdot [-1]-[a]\cdot [a]$.

\smallbreak

To show~(ii) it is by~(WK2) enough to show that $[b^{2}]=0$,
but this is a straightforward consequence of~(WK2), (WK4),
\ref{elemIdSubSect}~(2), and part~(i) of the corollary.
\end{proof}

Assertions~(i), (ii), and~(iv) of the following theorem are
less straightforward to prove. Note that~(ii) part~(b) and~(iv) are new
even for~$R$ a field.
\end{emptythm}

\begin{emptythm}
\label{KW-IdThm}
{\bfseries Theorem.}
{\it
\begin{itemize}
\item[(i)]
$\MW_{\ast}(R)$ is commutative.

\medbreak

\item[(ii)]
\begin{itemize}
\item[(a)]
$<a>_{\MW}+<b>_{\MW}=<a+b>_{\MW}+<ab(a+b)>_{\MW}$, and

\smallbreak

\item[(b)]
$[a]+[b]=[a+b]+[ab(a+b)]$
\end{itemize}

\noindent
for all $a,b\in R^{\times}$ such that also $a+b\in R^{\times}$.

\medbreak

\item[(iii)]
$[r]\cdot [st]+[s]\cdot [t]=[rs]\cdot [t]+[r]\cdot [s]$
for all $r,s,t\in R^{\times}$.

\smallbreak

\item[(iv)]
$[a+b]\cdot [ab(a+b)]\, =\, [a]\cdot [b]$
for all $a,b\in R^{\times}$ such that~$a+b\in R^{\times}$.
\end{itemize}
}

\begin{proof}
For~(i) it is enough to show that
$[a]\cdot [b]=[b]\cdot [a]$ for all $a,b\in R^{\times}$. By
Lemma~\ref{a-aRelLem} we have $[ab]\cdot [-ab]=0=[a]\cdot [-a]$ and so
using~(WK2) we get
$$
\begin{array}{r@{\, =\,}l}
0 & [ab]\cdot [-ab] \, =\, \big([a]+<a>_{\MW}\cdot [b]\big)\cdot
     \big([-a]+<-a>_{\MW}\cdot [b]\big) \\[2mm]
 & <a>_{\MW}\cdot\big([b]\cdot [-a]+<-1>_{\MW}\cdot [a]\cdot [b]\big)+
      <-a^{2}>_{\MW}\cdot [b]\cdot [b]\, .
\end{array}
$$
Using $[-a]=[a]+<a>_{\MW}\cdot [-1]$, $<-1>_{\MW}=-1$ and $<1>_{\MW}=1$
by~\ref{elemIdSubSect}~(1), and that $<rs^{2}>_{\MW}=<r>_{\MW}$ by
Corollary~\ref{a-1RelCor}~(ii) this equation is equivalent to
$$
0\, =\, <a>_{\MW}\cdot\big([b]\cdot [a]-[a]\cdot [b]\big)\, +\,
  [b]\cdot [-1]-[b]\cdot [b]\, .
$$
Now $[b]\cdot [b]=[b]\cdot [-1]$ by Corollary~\ref{a-1RelCor}~(i)
and so we get
$$
0\, =\, <a>_{\MW}\cdot\big([b]\cdot [a]-[a]\cdot [b]\big)\, .
$$ As $<a>_{\MW}$ is a unit in~$\MW_{\ast}(R)$ by~\ref{elemIdSubSect}~(1)
this proves our claim.

\smallbreak

Part~(a) of~(ii) can be proven as in Morel~\cite[Cor.\ 3.8]{Mo04a} for a field: By~(WK3)
we have $<r>_{\MW}+<1-r>_{\MW}=1+<r(1-r)>_{\MW}$ for~$r\in R^{\ast}$ with~$1-r$
also a unit. Setting $r=\frac{a}{a+b}$ gives then the result using that
$<rs^{2}>_{\MW}=<r>_{\MW}$ for all~$r,s\in R^{\times}$ by Corollary~\ref{a-1RelCor}~(ii).

\smallbreak

For part~(b) we set $c=a+b$. Then by~(WK2) and~(WK3), \ref{elemIdSubSect}~(3),
and Corollary~\ref{a-1RelCor} we get
$$
\begin{array}{r@{\, =\,}l}
   [abc\cdot c^{-1}] & [ab]\, =\, [abc^{-2}] \, =\, [ac^{-1}]+[bc^{-1}] \\[2mm]
     & [a]-<ac^{-1}>_{\MW}\cdot [c] +[b]-<bc^{-1}>_{\MW}\cdot [c]\, .
\end{array}
$$
On the other hand by~\ref{elemIdSubSect}~(3) again we have also
$$
[abc\cdot c^{-1}]\, =\, [abc]-<ab>_{\MW}\cdot [c]\, .
$$
Putting these two equations together we get using $<rs^{2}>_{\MW}=<r>_{\MW}$
that

\smallbreak

$\;\;\; [abc]+[c]-[a]-[b]$
\begin{equation}
\label{abcEq}
\begin{array}{c@{\;}l}
   = & <ab>_{\MW}\cdot [c]-\big(<a>_{\MW}+<b>_{\MW}\big)\cdot <c>_{\MW}\cdot [c]+[c] \\[2mm]
   = & <ab>_{\MW}\cdot [c]+\big(<a>_{\MW}+<b>_{\MW}\big)\cdot [c]+[c]\, ,
\end{array}
\end{equation}
where the second equality follows since  $<r>\cdot [r]=-[r]$ for all~$r\in R^{\times}$
as a direct computation using Corollary~\ref{a-1RelCor}~(i) shows.

\smallbreak

By part~(a) we have $<a>_{\MW}+<b>_{\MW}=<abc>_{\MW}+<c>_{\MW}$ and hence by~(\ref{abcEq})
$$
\begin{array}{r@{\; =\;}l}
  [abc]+[c]-[a]-[b]
      & <ab>_{\MW}\cdot [c]+\big(<abc>_{\MW}+<c>_{\MW}\big)\cdot [c]+[c] \\[2mm]
      & <ab>_{\MW}\cdot [c]+\big(<ab>_{\MW}+1\big)\cdot <c>_{\MW}\cdot [c] +[c] \\[2mm]
      & <ab>_{\MW}\cdot [c]-<ab>_{\MW}\cdot [c] -[c]+[c] \\[2mm]
      & 0\, ,
\end{array}
$$
as claimed.

\smallbreak

We show~(iii). This is a straightforward consequence of~(WK2) and the fact that
$\MW_{\ast}(R)$ is commutative:
$$
\begin{array}{r@{\, =\,}l}
  [rs]\cdot [t]+[r]\cdot [s] & \big( [s]+<s>_{\MW}\cdot [r]\big)\cdot [t] +[r]\cdot [s] \\[2mm]
   & [s]\cdot [t] +[r]\big( [s]+<s>_{\MW}\cdot [t]\big) \\[2mm]
   & [s]\cdot [t] +[r]\cdot [st]\, .
\end{array}
$$

\smallbreak

Finally we show~(iv). Setting in~(iii) $r=a+b$, $s=a(a+b)$,
and $t=b$, and using that $[xy^{2}]=[x]$ for all~$x,y\in R^{\times}$
by Corollary~\ref{a-1RelCor}~(ii) we get
\begin{equation}
\label{MW2Eq}
[a+b]\cdot [ab(a+b)]+[a(a+b)]\cdot [b]\, =\, [a]\cdot [b]+[a+b]\cdot [a(a+b)]\, .
\end{equation}
By~(WK2) and~(WK3) we have
$$
[x]\cdot [(1-x)y]\, =\, [x]\cdot [1-x]+[x]\cdot [y]+\eta\cdot [x]\cdot [1-x]\cdot [y]\, =\,
[x]\cdot [y]
$$
for all $x,y\in R^{\times}$, such that~$1-x$ also a unit in~$R$. This equation together with
Corollary~\ref{a-1RelCor}~(ii) and the fact that~$\MW_{\ast}(R)$ is commutative shows that
(note that $1-a(a+b)^{-1}=b(a+b)^{-1}$ is a unit in~$R$)
$$
\begin{array}{r@{\, =\,}l}
   [a(a+b)]\cdot [b] & [a(a+b)^{-1}]\cdot [b] \\[3mm]
      & [a(a+b)^{-1}]\cdot [(1-a(a+b)^{-1})b] \\[3mm]
      & [a(a+b)]\cdot [b^{2}(a+b)^{-1}]\, =\, [a+b]\cdot [a(a+b)]\, .
\end{array}
$$
Inserting this identity into~(\ref{MW2Eq}) proves the claimed identity.
We are done.
\end{proof}
\end{emptythm}

\begin{emptythm}
\label{KW-powersFdISubSect}
{\it Witt $K$-theory and the powers of the fundamental ideal.}
The {\it Witt algebra} of~$R$ is the $\Z$-graded $\W (R)$-algebra
$$
\WI_{\ast}(R)\, :=\;\bigoplus\limits_{n\in\Z}\WI_{n}(R)\, ,
$$
where $\WI_{n}(R)=\FdI^{n}(R)$ (recall that by convention
$\FdI^{n}(R)=\W(R)$ for~$n\leq 0$), with the obvious addition and
multiplication, \ie if~$x\in\WI_{m}(R)=\FdI^{m}(R)$ and~$y\in\WI_{n}(R)=\FdI^{n}(R)$
then $x\cdot y$ is the class of~$x\otimes y$ in~$\FdI^{m+n}(R)=\WI_{m+n}(R)$.
Following Morel~\cite{Mo04a} we set $\eta_{W}:=<1>\in\WI_{-1}(R)=\W(R)$. Then
multiplication by~$\eta_{W}$ corresponds to the natural embedding of~$\WI_{n+1}(R)$
into~$\WI_{n}(R)$ for all~$n\in\Z$.

\medbreak

We have $\ll u,1-u\gg=0$ for units~$u$ of~$R$, such that also
$1-u\in R^{\times}$,
$$
<1>\otimes\ll -1\gg\, =\, <1,1>\, =\, 2\;\in\ W (R)\, ,
$$
and $\ll ab\gg =\ll a\gg +\ll b\gg -\ll a,b\gg$ for all $a,b\in R^{\times}$.
Therefore there is a well defined homomorphism of $\Z$-graded $\Z$-algebras
$$
\Theta^{R}_{\ast}\, :\;\MW_{\ast}(R)\,\too\,\WI_{\ast}(R)\, ,\;
\begin{array}{r@{\,\,\longmapsto\,\,}l}
   [u] & \ll u\gg\,\in\WI_{1}(R)=\FdI (R) \\[3mm]
   \eta & \eta_{W}\,\in\WI_{-1}(R)=\W (R)\, .
\end{array}
$$
Since~$\FdI^{n}(R)$ for~$n\geq 1$ is additively generated by $n$-Pfister forms
and the $\Z$-module~$\W (R)$ is generated by $<a>=\Theta^{R}_{0}(<a>_{\MW})$, $a\in R^{\times}$,
we see that $\Theta^{R}_{n}$ is surjective for all~$n\geq 0$ and all local rings~$R$. Our
aim is to prove that~$\Theta^{R}_{\ast}$ is an isomorphism for all ``nice'' regular
local rings. More precisely, we show the following theorem, which is due
to Morel~\cite{Mo04a} if~$R$ is a field.
\end{emptythm}

\begin{emptythm}
\label{KW-mainThm}
{\bfseries Theorem.}
{\it
Let~$R$ be a field or a local ring which contains an infinite
field. Then
$$
\Theta^{R}_{\ast}\, :\;\MW_{\ast}(R)\,\too\,\WI_{\ast}(R)
$$
is an isomorphism.
}
\end{emptythm}

\goodbreak
\section{Proof of Theorem~\ref{KW-mainThm}}
\label{KW-mainThmProofSect}\bigbreak

\begin{emptythm}
\label{InStructureSubSect}
{\it A presentation of the powers of the fundamental ideal.}
Our proof of Theorem~\ref{KW-mainThm} uses the following presentation
of the powers of the fundamental ideal.

\medbreak

\noindent
{\bfseries Theorem.}
{\it
Let~$R$ be a field or a local ring whose residue field contains at least
$5$~elements. Let $\PfM_{n}(R)$ be the free abelian group generated by the
isometry classes of $n$-Pfister forms over~$R$. Following~\cite{Gi15} we
denote the isometry class of~$\ll a_{1},\ldots ,a_{n}\gg$ by $[a_{1},\ldots ,a_{n}]$.
\begin{itemize}
\item[(i)]
The kernel of the homomorphism
$$
\Z [R^{\times}]\,\too\,\W (R)\, ,\; [r]\,\longmapsto\, <r>
$$
is additively generated by $[1]-[-1]$, all $[ab^{2}]-[a]$ with $a,b\in R^{\times}$, and all expressions
$[a]+[b]-([a+b]+[ab(a+b)])$ with $a,b\in R^{\times}$, such that~$a+b\in R^{\times}$.

\smallbreak

\item[(ii)]
The kernel of the homomorphism
$$
\Z [R^{\times}]\,\too\,\FdI (R)\, ,\; [r]\,\longmapsto\, \ll r\gg
$$
is generated by $[1]$, all $[ab^{2}]-[a]$ with $a,b\in R^{\times}$, and
all expressions of the form $[a]+[b]-([a+b]+[ab(a+b)])$
with $a,b\in R^{\times}$, such that~$a+b\in R^{\times}$.

\smallbreak

\item[(iii)]
Assume now that~$R$ is a field or contains an infinite field of characteristic
not~$2$. Then for~$n\geq 2$ the kernel of the natural epimorphism
$$
\PfM_{n}(R)\,\too\,\FdI^{n}(R)\, ,\; [a_{1},\ldots ,a_{n}]\,\longmapsto\,
\ll a_{1},\ldots ,a_{n}\gg
$$
is generated by $[1,\ldots ,1]$; all sums
$$
[a,c_{2},\ldots ,c_{n}]+[b,c_{2},\ldots ,c_{n}]-\big( [a+b,c_{2},\ldots ,c_{n}]
+[ab(a+b),c_{2},\ldots ,c_{n}]\big)
$$
with $a,b\in R^{\times}$, such that $a+b\in R^{\times}$, and
$c_{2},\ldots ,c_{n}\in R^{\times}$;
and all sums
$$
[a,b,d_{3},\ldots ,d_{n}]+[ab,c,d_{3},\ldots ,d_{n}]-
\big( [b,c,d_{3},\ldots ,d_{n}]+[a,bc,d_{3},\ldots ,d_{n}]\big)
$$
with $a,b,c,d_{3},\ldots ,d_{n}\in R^{\times}$.
\end{itemize}
}

\begin{proof}
For~$R$ a local ring which contains an infinite field
this has been proven by the first author~\cite{Gi15}. For~$R$ a field
this is due to Witt~\cite{Wi37} if $n=0,1$ and to Arason and Elman~\cite[Thm.\ 3.1]{ArEl01}
if~$n\geq 2$. Note that the latter work as well as~\cite{Gi15} use the Milnor
conjectures which are now theorems by the work of Voevodsky~\cite{Voe03} and
Orlov, Vishik, and Voevodsky~\cite{OrViVoe07}.
\end{proof}

\smallbreak

\noindent
{\bfseries Remarks.}
\begin{itemize}
\item[(i)]
Morel~\cite[Lem.\ 3.10]{Mo04a} claims that if~$R$ is a field, then in the above presentation
of~$\FdI (R)$, the so-called Witt relation $[a]+[b]=[a+b]+[ab(a+b)]$ can be replaced
by the relation $[a]+[1-a]=[a(1-a)]$, $a\not=0,1$. As Detlev Hoffmann has
pointed out to us the Laurent field in one variable over $\Z/3\Z$ is a counter
example to this assertion.

\smallbreak

\item[(ii)]
Morel's~\cite{Mo04a} proof of Theorem~\ref{KW-mainThm} for a field~$R$
uses a slight alteration of Arason and Elman's~\cite[Thm.\ 3.1]{ArEl01}, the
assertion~\cite[Thm.\ 4.1]{Mo04a}. We do not know whether this follows
from Arason and Elman~\cite[Thm.\ 3.1]{ArEl01}. In any case, one can fixed this flaw
using our corollary to Theorem~\ref{n<=2Thm} below.
\end{itemize}

\medbreak

We start now the proof of Theorem~\ref{KW-mainThm}.
\end{emptythm}

\begin{emptythm}
\label{n<2CaseSubSect}
{\it The case~$n\leq 1$.}
Parts~(i) and~(ii) of the theorem above together with~\ref{elemIdSubSect}~(1),
Lemma~\ref{a-aRelLem}, Corollary~\ref{a-1RelCor}, and Theorem~\ref{KW-IdThm}~(i)
and~(ii) show that the morphisms
$$
\WI_{0}(R)\, =\,\W (R)\,\too\,\MW_{0}(R)\, ,\;\; <a>\,\longmapsto\, <a>_{\MW}
$$
and
$$
\WI_{1}(R)\, =\,\FdI (R)\,\too\,\MW_{1}(R)\, ,\;\;\ll a\gg\,\longmapsto\, [a]
$$
are well defined. These are obviously inverse to~$\Theta^{R}_{0}$
and~$\Theta^{R}_{1}$, respectively.

\smallbreak

It follows from this that~$\Theta^{R}_{n}$ is also an isomorphism for~$n<0$.
In fact, we have a commutative diagram
$$
\xymatrix{
{\MW_{n+1}(R)} \ar[r]^-{\cdot\eta} \ar[d]_-{\Theta^{R}_{n+1}} & {\MW_{n}(R)}
   \ar[d]^-{\Theta^{R}_{n}}
\\
{\WI_{n+1}(R)} \ar[r]_-{\cdot\eta_{W}}^-{\simeq} & {\WI_{n}(R)}
}
$$
for all $n<0$. By induction we can assume that $\Theta^{R}_{n+1}$ is an
isomorphism. Since $\MW_{n+1}(R)\xrightarrow{\cdot\eta}\MW_{n}(R)$ is an
epimorphism for all~$n<0$ this implies that also $\Theta^{R}_{n}$ is an
isomorphism.
\end{emptythm}

\begin{emptythm}
\label{n=2CaseSubSect}
{\it The case~$n=2$.}
The following argument is an adaption of a trick of
Suslin~\cite[Proof of Lem.\ 6.3~(iii)]{Su87} to Witt $K$-theory.

\smallbreak

Consider the product of sets
$$
\MW_{2}(R)\times\; R^{\times}/(R^{\times})^{2}\, .
$$
We define an addition on this set by the rule
$$
(x,\bar{r})+(y,\bar{s})\, :=\; (x+y+[r]\cdot [s],\bar{r}\cdot\bar{s})\, ,
$$
where~$\bar{r}$ denotes the class of~$r\in R^{\times}$
in~$R^{\times}/(R^{\times})^{2}$. As $[ab^{2}]=[a]$ in~$\MW_{1}(R)$
by Corollary~\ref{a-1RelCor}~(ii) this is well defined. It is straightforward
to check that with this addition the product of sets
$\MW_{2}(R)\times\; R^{\times}/(R^{\times})^{2}$ is an abelian group with~$(0,\bar{1})$
as zero and inverse $-(x,\bar{r})=(-x-[r]\cdot [r^{-1}],\bar{r}^{-1})$.

\smallbreak

Using part~(ii) of the theorem in~\ref{InStructureSubSect} together
with~\ref{elemIdSubSect}~(1), Corollary~\ref{a-1RelCor}~(ii), and
Theorem~\ref{KW-IdThm}~(iv) we see that
$$
\FdI (R)\,\too\,\MW_{2}(R)\times\; R^{\times}/(R^{\times})^{2}\, ,\;\;
\ll r\gg\,\longmapsto\,\, (0,\bar{r})
$$
is a well defined homomorphism. The image of the restriction of this map
to~$\FdI^{2}(R)\subseteq\FdI (R)$ is contained in the subgroup
$$
\big\{\, (x,\bar{1})\, |\, x\in\MW_{2}(R)\,\big\}
$$
of $\MW_{2}(R)\times\; R^{\times}/(R^{\times})^{2}$ which is naturally
isomorphic to~$\MW_{2}(R)$ via $(x,\bar{1})\mapsto x$. The induced
homomorphism $\FdI^{2}(R)\too\MW_{2}(R)$ maps $\ll r,s\gg$ to~$[r]\cdot [s]$
and is therefore inverse to $\Theta^{R}_{2}:\MW_{2}(R)\too\WI_{2}(R)=\FdI^{2}(R)$.

\smallbreak

Altogether we have shown:
\end{emptythm}

\begin{emptythm}
\label{n<=2Thm}
{\bfseries Theorem.}
{\it
Let~$R$ be a field or a local ring whose residue field has at least $5$~elements.
Then
$$
\Theta_{n}^{R}\, :\;\MW_{n}(R)\,\too\,\WI_{n}(R)\, =\,\FdI^{n}(R)
$$
is surjective for all $n\in\Z$ and an isomorphism for all $n\leq 2$.
}

\medbreak

\noindent
{\bfseries Corollary.}
{\it
Let~$R$ be a field or a local ring whose residue field has at least $5$~elements,
and~$n\geq 1$. Then
$$
\ll a_{1},\ldots ,a_{n}\gg\,\simeq\,\ll b_{1},\ldots ,b_{n}\gg
$$
if and only if
$$
[a_{1}]\cdot\ldots\cdot [a_{n}]\, =\, [b_{1}]\cdot\ldots\cdot [b_{n}]
$$
in~$\MW_{n}(R)$ for all $a_{1},\ldots ,a_{n},b_{1},\ldots ,b_{n}\in R^{\times}$.
In particular, $[a_{1}]\cdot\ldots [a_{n}]=0$ in~$\MW_{n}(R)$ if and only if
$\ll a_{1},\ldots ,a_{n}\gg=0$ in~$\W (R)$.
}

\begin{proof}
If~$n=1$ this is obvious by the theorem above, so let~$n\geq 2$. By
Theorem~\ref{p-chainThm} we know that isometric $n$-Pfister forms are chain
$p$-equivalent and so we are reduced to the case $n=2$, where the assertion
follows again from the theorem above.
\end{proof}

\smallbreak

With this corollary we can finish the proof of Theorem~\ref{KW-mainThm}.
\end{emptythm}

\begin{emptythm}
\label{endProofSubSect}
{\it End of the proof of Theorem~\ref{KW-mainThm}.}
Let now~$n\geq 3$ and~$R$ be a field or a local
ring which contains an infinite field. It follows from the
corollary to Theorem~\ref{n<=2Thm} above that the homomorphism
$$
\PfM_{n}(R)\,\too\,\MW_{n}(R)\, ,\;\;
[a_{1},\ldots ,a_{n}]\,\longmapsto\, [a_{1}]\cdot\ldots\cdot [a_{n}]
$$
is well defined. Part~(iii) of the theorem in~\ref{InStructureSubSect}
together with~\ref{elemIdSubSect}~(1) and Theorem~\ref{KW-IdThm}~(ii) and~(iii)
shows that this homomorphism factors through~$\FdI^{n}(R)=\WI_{n}(R)$.
It is obviously inverse to~$\Theta^{R}_{n}$. We are done.
\end{emptythm}

\goodbreak
\section{Milnor-Witt $K$-theory of a local ring}
\label{MWKSect}\bigbreak

\begin{emptythm}
\label{MWKDef}
{\bfseries Definition.}
Let~$R$ be a local ring. The {\it Milnor-Witt $K$-ring of~$R$} is quotient of the graded
and free $\Z$-algebra generated by elements~$\{ a\}$, $a\in R^{\times}$, in degree~$1$
and one element~$\hat{\eta}$ in degree~$-1$ by the two sided ideal which is generated by
the expressions
\begin{itemize}
\item[(MW1)]
$\hat{\eta}\cdot\{ a\}-\{ a\}\cdot\hat{\eta}$, $a\in R^{\times}$;

\smallbreak

\item[(MW2)]
$\{ ab\}-\{ a\}-\{ b\} +\hat{\eta}\cdot\{ a\}\cdot\{ b\}$, $a,b\in R^{\times}$;

\smallbreak

\item[(MW3)]
$\{ a\}\cdot\{ 1-a\}$, if $a$ and~$1-a$ in~$R^{\times}$; and

\smallbreak

\item[(MW4)]
$\hat{\eta}\cdot\big( 2+\hat{\eta}\cdot\{ -1\}\big)$.
\end{itemize}
We denote this graded $\Z$-algebra by $\MWK_{\ast}(R)=\bigoplus\limits_{n\in\Z}\MWK_{n}(R)$.
Note that by~(MW2) for~$n\geq 1$ the $n$th graded piece~$\MWK_{n}(R)$ is generated as abelian
group by all products $\{ a_{1},\ldots ,a_{n}\}\:=\{ a_{1}\}\cdot\ldots\cdot\{ a_{n}\}$.

\medbreak

We prove now our main result about Milnor-Witt $K$-theory of local rings
following Morel's arguments~\cite[Sect.\ 5]{Mo04a} in the field case.
\end{emptythm}

\begin{emptythm}
\label{WK-MK-MWKSubSect}
{\it Milnor $K$-theory, Witt $K$-theory, and Milnor-Witt $K$-theory.}
Let~$R$ be a local ring. Milnor $K$-theory of a field
has been introduced by Milnor~\cite{Mi69/70}. The definition for rings
we use here is the naive generalization of the one by Milnor and seems
to have appeared for the first time in Nesterenko and Suslin~\cite{NeSu90}.

\smallbreak

Denote by~$\tensor_{\Z}\,(R^{\times})$ the tensor algebra of the abelian
group~$R^{\times}$ over~$\Z$. The {\it Milnor $K$-theory} of~$R$ is the quotient
of this algebra by the ideal generated by all tensors $a\otimes (1-a)$ with~$a$
and~$1-a$ units in~$R$. This is a graded $\Z$-algebra which we denote
$\MK_{\ast}(R)=\bigoplus\limits_{n\geq 0}\MK_{n}(R)$.
The class of a tensor $a_{1}\otimes\ldots\otimes a_{n}$ will be denoted
as in the original source Milnor~\cite{Mi69/70} by~$\ell (a_{1})\cdot\ldots\cdot\ell (a_{n})$
as we have the now usual $\{ a_{1},\ldots ,a_{n}\}$ reserved for symbols in Milnor-Witt
$K$-theory.

\smallbreak

There is a natural surjective homomorphism of graded $\Z$-algebras
$$
\varpi_{\ast}^{R}\, :\;\MWK_{\ast}(R)\,\too\,\MK_{\ast}(R)\, ,
$$
which maps~$\hat{\eta}$ to~$0$ and~$\{ a\}$
to~$\ell (a)$, $a\in R^{\times}$. The ideal
$\hat{\eta}\cdot\MWK_{\ast +1}(R)$ is
in the kernel, and the induced homomorphism
$\MWK_{n}(R)/\hat{\eta}\cdot\MWK_{n+1}(R)\too\MK_{\ast}(R)$
is an isomorphism. The inverse maps the symbol
$\ell (a_{1})\cdot\ldots\cdot\ell (a_{n})$ to
$\{ a_{1}\}\cdot\ldots\cdot\{ a_{n}\}$ modulo~$\hat{\eta}\cdot\MWK_{n+1}(R)$.

\smallbreak

Set now $h:=2+\hat{\eta}\cdot\{ -1\}\in\MWK_{0}(R)$. Then by~(MW4) we
have $\hat{\eta}\cdot h=0$ and so there is an exact sequence
$$
\MWK_{n+1}(R)/h\cdot\MWK_{n+1}(R)\,\xrightarrow{\;\cdot\hat{\eta}\;}\,\MWK_{n}(R)\,
\xrightarrow{\;\cdot\varpi_{n}^{R}\;}\,\MK_{n}(R)\,\too\, 0
$$
for all~$n\in\Z$, where $\MK_{n}(R)=0$ for~$n<0$ is understood. The group
on the left hand side of this sequence is isomorphic to~$\MW_{n}(R)$. In fact,
$\hat{\eta}\mapsto\eta$ and $\{ u\}\mapsto -[u]$, $u\in R^{\times}$, defines
a morphism of $\Z$-graded $\Z$-algebras $\MWK_{\ast}(R)\too\MW_{\ast}(R)$, whose
kernel contains the ideal generated by~$h\in\MWK_{0}(R)$, \ie we have a homomorphism
of $\Z$-graded $\Z$-algebras
$\MWK_{\ast}(R)/h\cdot \MWK_{\ast}(R)\too\MW_{\ast}(R)$.
This is an isomorphism. The inverse maps~$[u]$ to~$-\{ u\}$
and~$\eta$ to~$\hat{\eta}$ modulo $h\cdot\MWK_{\ast}(R)$. Hence we have
an exact sequence
\begin{equation}
\label{MWKExSeq}
\MW_{n+1}(R)\,\xrightarrow{\;\epsilon^{R}_{n}\;}\,\MWK_{n}(R)\,
\xrightarrow{\;\varpi^{R}_{n}\;}\,\MK_{n}(R)\,\too 0
\end{equation}
for all~$n\in\Z$, where $\epsilon^{R}_{n}$ maps
$\eta^{r}\cdot [u_{1}]\cdot\ldots\cdot [u_{n+r}]$
to $(-1)^{n+r+1}\hat{\eta}^{r+1}\cdot\{ u_{1},\ldots ,u_{n+r}\}$.
\end{emptythm}

\begin{emptythm}
\label{MWK-FdISubSect}
{\it Milnor-Witt $K$-theory and the powers of the fundamental ideal.}
Let~$R$ be a local ring. As for Witt $K$-theory there is also a homogeneous
homomorphism of $\Z$-graded $\Z$-algebras
$$
\Upsilon^{R}_{\ast}\, :\;\MWK_{\ast}(R)\,\too\,\WI_{\ast}(R)\, ,\;\;
\begin{array}{r@{\,\,\longmapsto\,\,}l}
   \{ u\} & -\ll u\gg\,\in\WI_{1}(R)=\FdI (R) \\[3mm]
   \hat{\eta} & \eta_{W}\,\in\WI_{-1}(R)=\W (R)\, .
\end{array}
$$
We leave it to the reader to check that this map is well defined. It is
also surjective as~$\Theta^{R}_{\ast}:\MW_{\ast}(R)\too\WI_{\ast}(R)$
and~$\Upsilon^{R}_{\ast}$ have the same image.

\smallbreak

These maps fit into the following commutative diagram with exact rows
\begin{equation}
\label{cartSqDiag}
\xymatrix{
 & \MW_{n+1}(R) \ar[r]^-{\epsilon^{R}_{n}} \ar[d]_-{\Theta^{R}_{n+1}} & \MWK_{n}(R)
   \ar[r]^-{\varpi^{R}_{n}} \ar[d]_-{\Upsilon^{R}_{n}} & \MK_{n}(R) \ar[r]
       \ar[d]^-{e^{R}_{n}} & 0
\\
0 \ar[r] & \FdI^{n+1} (R) \ar[r]_-{\subseteq} & \FdI^{n} (R) \ar[r] &
   \FdI^{n}(R)/\FdI^{n+1}(R) \ar[r] & 0 \rlap{\, ,}
}
\end{equation}
where
$$
e_{n}^{R}:\MK_{n}(R)\,\too\,\FdI^{n}(R)/\FdI^{n+1}(R)\, ,\;\;
\ell (a_{1})\cdot\ldots\cdot\ell (a_{n})\,\longmapsto\,\ll a_{1},\ldots, a_{n}\gg +\FdI^{n+1}(R)\, .
$$
By Theorems~\ref{KW-mainThm} and~\ref{n<=2Thm} we conclude by a
diagram chase on~(\ref{cartSqDiag}) the following result.
\end{emptythm}

\begin{emptythm}
\label{MWK-mainThm}
{\bfseries Theorem.}
{\it
Let~$R$ be a field or a local ring whose residue field has at least $5$~elements.
Then
$$
\xymatrix{
\MWK_{n}(R) \ar[r]^-{\varpi^{R}_{n}} \ar[d]_-{\Upsilon^{R}_{n}} &
    \MK_{n}(R) \ar[d]^-{e^{R}_{n}}
\\
\FdI^{n}(R) \ar[r] & \FdI^{n}(R)/\FdI^{n+1}(R)
}
$$
is a pull-back diagram in the following cases:
\begin{itemize}
\item[(a)]
$n\leq 1$; or

\smallbreak

\item[(b)]
$n$~is arbitrary and~$R$ is a field or a local ring which contains an
infinite field.
\end{itemize}
In particular, we have natural isomorphisms~$\MWK_{0}(R)\simeq\GW (R)$
(by the lemma in~\ref{GW-Z/2LemSubSect}) and $\MWK_{n}(R)\simeq\W (R)$ for all $n<0$,
and if~$R$ is a field or a regular local ring which contains an infinite
field then the sequence
$$
\xymatrix{
0 \ar[r] & \MW_{n+1}(R) \ar[r]^-{\epsilon^{R}_{n}} & \MWK_{n}(R)
   \ar[r]^-{\varpi^{R}_{n}} & \MK_{n}(R) \ar[r] & 0
}
$$
is exact for all~$n\in\Z$.
}
\end{emptythm}

\goodbreak
\section{Unramified Milnor Witt $K$-groups}
\label{UnrMWKSect}\bigbreak

\begin{emptythm}
\label{residueMapsSubSect}
{\it Residue maps.}
Let~$F$ be a field with discrete valuation~$\nu$ and residue field~$F(\nu)$.
Denote by~$\pi$ a uniformizer for~$\nu$. As well known, there is a so called
(second) residue homomorphism
$$
\partial^{\nu,\pi}_{W}\, :\;\W (F)\,\too\,\W (F(\nu))\, ,\;
<u\pi^{i}>\,\longmapsto\, i\cdot <\bar{u}>\, ,
$$
where~$u$ is a $\nu$-unit, $i\in\{ 0,1\}$, and~$\bar{u}$ denotes the image of~$u$
in the residue field~$F(\nu)$. This homomorphism depends on the choice of the
uniformizer~$\pi$, but its kernel is independent of this choice. It has
been shown by Milnor~\cite[\S 5]{Mi69/70} that
$\partial^{\nu,\pi}_{W}(\FdI^{n}(F)\subseteq\FdI^{n-1}(F)$ for all~$n\in\N$.
Let $\grFdI_{n}(R):=\FdI^{n}(R)/\FdI^{n+1}(R)$. Then$~\partial^{\nu,\pi}_{W}$
induces a residue homomorphism
$$
\partial^{\nu}_{\grFdI}\, :\;\grFdI_{n}(F)\,\too\,\grFdI_{n-1}(F(\nu))
$$
which does not depend on the choice of~$\pi$.

\smallbreak

Similarly there is a (second) residue homomorphism~$\partial^{\nu}_{M}$:
$$
\MK_{n}(F)\,\too\,\MK_{n-1}(F(\nu))\, ,\;
\ell (u\pi^{i})\cdot\ell (u_{2})\cdot\ldots\cdot\ell (u_{n})\,\longmapsto\,
i\cdot\ell (\bar{u}_{2})\cdot\ldots\cdot\ell(\bar{u}_{n})\, ,
$$
where $u,u_{2},\ldots ,u_{n}$ are $\nu$-units and $i\in\{ 0,1\}$. This map
does not depend on~$\pi$.

\smallbreak

Since~$\MWK_{n}(F)$ is the pull-back of~$\FdI^{n}(F)$ and~$\MK_{n}(F)$
over~$\grFdI_{n}(F)$ by Theorem~\ref{MWK-mainThm}, these residue maps
induce a {\it residue homomorphism} on Milnor-Witt $K$-theory
$$
\partial^{\nu,\pi}_{MW}\, :\;\MWK_{n}(F)\,\too\,\MWK_{n-1}(F(\nu))\, ,\;
$$
which is $\eta_{MW}$-linear and is uniquely determined by
$$
\partial^{\nu,\pi}_{MW}(\{ \pi ,u_{2},\ldots ,u_{n}\})\, =\,
\{ \bar{u}_{2},\ldots ,\bar{u}_{n}\}
$$
and $\partial^{\nu,\pi}_{MW}(\{ u_{1},u_{2},\ldots ,u_{n}\})=0$,
where $u_{1},\ldots ,u_{n}$ are $\nu$-units. Hence it coincides with the residue map
constructed by Morel~\cite[Thm.\ 3.15]{A1AlgTop}. It depends on the choice of~$\pi$,
but the kernel does not.
\end{emptythm}

\begin{emptythm}
\label{unrMWKSubSect}
{\it Unramified Milnor-Witt $K$-groups.}
Let~$X$ be an integral  locally noetherian scheme which is regular in codimension one.
Then~$\OXx$ is a discrete valuation ring for all $x\in X^{(1)}$, the set
of points of codimension one in~$X$. Choose an uniformizer~$\pi_{x}\in\OXx$
for all such~$x$. Then we have residue maps
$$
\partial^{\pi_{x}}_{MW}\, :\;\MWK_{n}(F(X))\,\too\,\MWK_{n-1}(F(x))\, ,
$$
where~$F(X)$ is the function field of~$X$ and~$F(x)$ the residue field of~$x\in X^{(1)}$.
The $n$th {\it unramified Milnor-Witt $K$-group} of~$X$ is defined as
$$
\MWK_{n,unr}(X)\, :=\;\Ker\Big(\,\MWK_{n}(F(X))\,
\xrightarrow{\;\sum_{x\in X^{(1)}}\partial^{\pi_{x}}_{MW}\;}\bigoplus\limits_{x\in X^{(1)}}\MWK_{n-1}(F(x))\,\Big)\, .
$$
The unramified groups~$\Wunr (X)$, $\FdIunr^{n}(X)$, $\grFdI_{n,unr}(X)$, and~$\MK_{n,unr}(X)$
are defined analogously. If $X=\Spec R$ is affine we use $\MWK_{n,unr}(R)$ instead
of~$\MWK_{n,unr}(X)$, $\Wunr (R)$ instead of~$\Wunr (X)$ and so on.

\smallbreak

Assume now that~$X$ is the spectrum of a regular local ring~$R$ which contains an
infinite field. Denote by~$K$ the quotient field of~$R$. It has been shown
in~\cite{BaGiPaWa02} that $\W (R)\too\W(K)$ induces an isomorphism
$\W(R)\xrightarrow{\simeq}\Wunr (R)$. Kerz and M\"uller-Stach~\cite[Cor.\ 0.5]{KeMS07}
have shown $\FdI^{n}(K)\cap\Wunr (R)=\FdI^{n}(R)$ for all integers~$n$ and therefore
we have also natural isomorphisms $\FdI^{n}(R)\xrightarrow{\simeq}\FdIunr^{n}(R)$
and $\grFdI_{n}(R)\xrightarrow{\simeq}\grFdI_{n,unr}(R)$.
On the other hand, the main result of Kerz~\cite{Ke09} asserts that
$\MK_{n}(R)\too\MK_{n,unr}(R)$ is also an isomorphism for all~$n\in\N$.
These isomorphisms together with our Theorem~\ref{MWK-mainThm} imply the
following result.
\end{emptythm}

\begin{emptythm}
\label{unr-mainThm}
{\bfseries Theorem.}
{\it
Let~$R$ be a regular local ring which contains an infinite field.
Then the natural homomorphism $\MWK_{n}(R)\too\MWK_{n,unr}(R)$ is an isomorphism, \ie
the sequence
$$
\xymatrix{
0 \ar[r] & \MWK_{n}(R) \ar[r] & \MWK_{n}(K) \ar[rr]^-{\sum_{\hgt P=1}\partial^{\pi_{P}}_{MW}} & &
    {\bigoplus\limits_{\hgt P=1}} \MWK_{n-1}(R_{P}/PR_{P})\, ,
}
$$
where~$K$ is the quotient field of~$R$, is exact for all integers~$n$.
}

\smallbreak

As observed by Colliot-Th\'el\`ene~\cite{CT95} this has the following consequence.
\end{emptythm}

\begin{emptythm}
\label{birMWCor}
{\bfseries Corollary.}
{\it
The $n$th unramified Milnor Witt $K$-group is a birational invariant
of smooth and proper $F$-schemes for all integers~$n$ and infinite
fields~$F$.
}

\begin{proof}
For the convenience of our reader we recall briefly Colliot-Th\'el\`ene's
argument. Let $f:X\dashrightarrow Y$ be a birational morphism between
smooth and proper $F$-schemes. By symmetry it is enough to show that
the induced homomorphism $f^{\ast}:\MWK_{n}(F(Y))\too\MWK_{n}(F(X))$ maps
unramified elements to unramified elements. To see this we observe
first that since~$Y$ is proper we can assume that~$f$ is defined on
an open $U\subseteq X$ which contains~$X^{(1)}$. Let now~$x\in X^{(1)}$
and~$y=f(x)$. Then we have a commutative diagram
$$
\xymatrix{
\MWK_{n}(F(Y)) \ar[r]^-{f^{\ast}} & \MWK_{n}(F(X))
\\
\MWK_{n}(\OYy) \ar[u] \ar[r]_-{f^{\ast}} &
\MWK_{n}(\OXx) \ar[u]
}
$$
from which the claim follows by the theorem above.
\end{proof}
\end{emptythm}

\begin{emptythm}
\label{WittGroupRemSubSect}
{\bfseries A remark on Witt groups of schemes.}
The same method applied to Witt groups shows that the unramified
Witt group of a smooth and proper scheme over a field of characteristic not~$2$
is also a birational invariant. This implication of the Gersten conjecture for
Witt groups was already observed by Colliot-Th\'el\`ene in~\cite{CT95} but
at that time the Gersten conjecture was only known for low dimensional
regular local rings. Only later Balmer~\cite{Ba01} proved it for regular
local rings of geometric type and Balmer, Panin, Walter and the first author~\cite{BaGiPaWa02}
proved it for regular local rings containing a field (of characteristic not~$2$).

\smallbreak

We finally mention the following application of the fact that the unramified
Witt group is a birational invariant of smooth and proper schemes over a
field of characteristic~$\not= 2$, which follows from this by the
Balmer-Walter~\cite{BaWa02} spectral sequence.

\medbreak

\noindent
{\bfseries Theorem.}
{\it
The Witt group is a birational invariant of smooth and proper schemes of
dimension~$\leq 3$ over a field of characteristic not~$2$.
}

\medbreak

\noindent
By a theorem of Arason~\cite{Ar80} we known that~$W(F)\simeq\W (\P^{n}_{F})$
for all positive integers~$n$, and so
$$
\W(F)\,\simeq\W (X)
$$
for all rational smooth and proper $F$-schemes~$X$ of dimension~$\leq 3$.
\end{emptythm}

\bibliographystyle{amsalpha}

\end{document}